\documentclass[a4paper, 11pt, preprint, 3p, times]{P-article}
\usepackage{amssymb}
\usepackage{epsfig}
\usepackage{amsfonts}
\usepackage{amsmath, array}
\usepackage{euscript}
\usepackage{amscd}
\usepackage{amsthm}
\usepackage{ulem}
\usepackage{xcolor}
\usepackage{mathrsfs}
\usepackage{eufrak}
\usepackage{yfonts}
\usepackage{booktabs}

\usepackage[english]{babel}
\usepackage{blindtext}
\usepackage{mathtools}
\usepackage{xcolor}
\usepackage{soul} 

\usepackage[pdftex,
pdfauthor=Prosenjit Das,
pdfsubject={Affine Algebraic Geometry, Commutative Algebra},
pdfproducer=TeXStudio,
pdfcreator=pdflatex]{hyperref}

\usepackage{hyperref}

\hypersetup{
	colorlinks=true,
}

\DeclareMathAlphabet{\mathpzc}{OT1}{pzc}{m}{it}
\newtheorem{thm}{Theorem}[section]
\newtheorem{lem}[thm]{Lemma}
\newtheorem{prop}[thm]{Proposition}

\newtheorem{cor}[thm]{Corollary}

\newdefinition{defn}[thm]{Definition}
\newdefinition{ex}[thm]{Example}
\newdefinition{rem}[thm]{Remark}
\newdefinition{note}{Note}

\newcommand{\comment}[1]{}

\newcommand\m {\mathfrak{m}}

\newcommand{\A}[1]{\mathbb{A}^{#1}}
\newcommand{\KerD}{\text{Ker}(D)}

\newcommand{\Ker}[1]{\text{Ker}(#1)}

\newcommand{\Sym}[2]{\text{Sym}_{#1}(#2)}

\newcommand{\Ht}[1]{\text{ht}(#1)}
\newcommand{\ol}[1] {\overline{#1}}
\newcommand{\UL}[1] {\underline{#1}}

\newcommand{\Qt}[1]{\text{Qt}(#1)}

\newcommand{\Spec}[1]{\text{Spec}(#1)}
\newcommand{\mbbQ}{\mathbb{Q}}

\newcommand{\mbbN}{\mathbb{N}}
\newcommand{\trdeg}[2]{\text{tr.deg}_{#1}(#2)}

\newcommand{\JacDer}[2]{\mathcal{J}\mathcal{D}_{(#2)}(#1, -)}
\newcommand{\JacMat}[2]{\mathcal{J}_{{(#2)}}(#1)}

\begin{document}
	\begin{frontmatter}
		\title{A criterion to determine residual coordinates of $\mathbb{A}^2$-fibrations}
		
		\author{Janaki Raman Babu}
		\address{Department of Mathematics, Indian Institute of Space Science and Technology, \\
			Valiamala P.O., Trivandrum 695 547, India\\
			email: \texttt{raman.janaki93@gmail.com, janakiramanb.16@res.iist.ac.in}}
		
		\author{Prosenjit Das\footnote{Corresponding author.}}
		\address{Department of Mathematics, Indian Institute of Space Science and Technology, \\
			Valiamala P.O., Trivandrum 695 547, India\\
			email: \texttt{prosenjit.das@gmail.com, prosenjit.das@iist.ac.in}}
		
			\begin{abstract}
			This article discusses a criterion to determine residual coordinates of an $\mathbb{A}^2$-fibration over a Noetherian domain containing $\mathbb{Q}$.
		\end{abstract}
		\end{frontmatter}

	\section{Introduction} \label{Sec_Intro}
	Throughout this article rings will be commutative with unity. Let $R$ be a ring. The \textit{polynomial ring} in $n$ variables over $R$ is denoted by $R^{[n]}$. Let $A$ be an $R$-algebra. We shall use the notation $A = R^{[n]}$ to mean that $A$ is isomorphic, as an $R$-algebra, to a polynomial ring in $n$ variables over $R$. $A$ is called a \textit{stably polynomial} algebra over $R$, if $A^{[m]} = R^{[n]}$ for some $m,n \in \mbbN$. A tuple $(X_1, X_2, \cdots, X_r)$ of elements of $A = R^{[n]}$ is called a \textit{coordinate-tuple} of $A$ over $R$ if $A = R[X_1, X_2, \cdots, X_r]^{[n-r]}$. For a prime ideal $P$ of $R$, let $k(P)$ denote the \textit{residue field} $R_P/PR_P$. $A$ is called an \textit{$\A{n}$-fibration} or \textit{affine $n$-fibration} over $R$, if $A$ is finitely generated and flat over $R$, and $A \otimes_R k(P) = k(P)^{[n]}$ for all $P \in \Spec{R}$. A tuple of algebraically independent elements $(W_1, W_2, \cdots, W_r)$ of an $\A{n}$-fibration $A$ over $R$ is called a \textit{residual coordinate tuple} of $A$ if $A \otimes_R k(P) = (R[W_1, W_2, \cdots, W_r] \otimes_R k(P))^{[n-r]}$ for each $P \in \Spec{R}$. Let $A = R^{[n]}$ and $(X_1, X_2 \cdots, X_n)$ is a coordinate system of $A$, then $\JacMat{F_1, F_2, \cdots, F_{m}}{X_1, X_2, \cdots, X_n}$ shall denote the \textit{Jacobian matrix} of $F_1, F_2, \cdots, F_{m} \in A$ with respect to the coordinate system $X_1, X_2, \cdots, X_n$ and $\JacDer{F_1, F_2, \cdots, F_{n-1}}{X_1, X_2, \cdots, X_n}$ shall denote the \textit{Jacobian derivation} induced by $F_1, F_2, \cdots, F_{n-1} \in A$ with respect to the coordinate system $X_1, X_2, \cdots, X_n$. 
	
	\medskip
	
	Given an element of a polynomial algebra, checking whether it is a coordinate is a difficult problem, especially when the base ring not a field. The problem is not easy even for the two-variable polynomial algebras. Since any coordinate is necessarily a residual coordinate, to determine whether an element is a coordinate one should at least check if it is a residual coordinate. The theory of residual coordinates of polynomial algebras, which is introduced by Bhatwadekar-Dutta in \cite{BD_RES}, asserts that if the base ring is Noetherian containing $\mbbQ$, then any residual coordinate of a polynomial algebra in two variables is a coordinate (see \cite[Theorem 3.2]{BD_RES}). However, determining whether an element is a residual coordinate requires a significant effort. Bhatwadekar-Dutta, in \cite{BD_LND}, established the following useful criterion for an element of a polynomial algebra in two variables over a Noetherian domain containing $\mathbb{Q}$ to be a residual coordinate, and hence a coordinate.
	
	\begin{thm} \label{BD_variable-check}
	Let $R$ be a Noetherian domain containing $\mbbQ$ with quotient field $K$ and $F \in R[X,Y] = R^{[2]}$ be such that $K[X,Y] = K[F]^{[1]}$. Then, the following are equivalent.
	
	\begin{enumerate}
		\item [\rm (I)] $(F_X, F_Y)R[X,Y] = R[X,Y]$.
		\item [\rm (II)] $F$ is a residual coordinate of $R[X,Y]$.
		\item [\rm (III)] $F$ is a coordinate of $R[X,Y]$, i.e., $R[X,Y] = R[F]^{[1]}$.
	\end{enumerate}
	\end{thm}  
	
	Since the problem of affine fibration asks whether a given affine fibration is a polynomial algebra, it remained as one of the motivating forces to establish a theory of residual coordinates of affine fibrations with the aim to judge whether an element of an affine fibration is a coordinate. Kahoui-Ouali in \cite{Kahoui_Residual} and Das-Dutta in \cite{DD_residual} developed the concept and the theory of residual coordinates of affine fibrations which subsequently helped as a tool to tackle many problems in affine fibrations, e.g., \cite{Das_cancel}, \cite{Kahoui_A2-fib_triviality-criterion}, \cite{Lahiri_partial-coordinate-system}, \cite{Kahoui-Essamaoui-Ouali_ResCoord-one-dim}, \cite{Babu-Das_Struct_A2-fib_FPF-LND}.  
	\medskip
	
	In Section \ref{Sec_Recognizing-residual-variables} of this article we establish the following sufficient criterion (see Theorem \ref{Thm_Recognize-ResCord-of-A2Fib-RingCase}) which help determine whether an element of an $\A{2}$-fibration is a residual coordinate. The result is an analogue of Theorem \ref{BD_variable-check} in the context of $\A{2}$-fibration.
	
	\smallskip
	
	\noindent
	\textbf{Theorem A.} \textit{Let $R$ be a Noetherian ring containing $\mathbb{Q}$ and $A$ an $R$-algebra such that $A$ is a retraction of $B = R[X_1, X_2, \cdots, X_{n}] = R^{[n]}$ and $\trdeg{R/P}{A \otimes_{R} R/P} =2$ for each minimal prime ideal $P$ of $R$. Let $F \in A$ be such that $A\otimes_R \Qt{R/P} = \Qt{R/P}[F]^{[1]}$ for each for each minimal prime ideal $P$ of $R$. If $(F_{X_1},F_{X_2}, \cdots, F_{X_{n}})B = B$, then $F$ is a residual coordinate of $A$.}
	
	\medskip
	
	One must note that the $R$-algebra $A$ in Theorem A is an $\A{2}$-fibration over $R$ due to a result of Chakraborty-Dasgupta-Dutta-Gupta in \cite{Neena-Sagnik-Nikhilesh_Retract-of-polynomial} (see Corollary \ref{Cor_Neena-Nikhilesh-Sagnik_Retraction-A2Fib}).  As a corollary of Theorem A we get a necessary and sufficient condition for an element of a stably polynomial $\A{2}$-fibration to be a residual coordinate. 
	
	\smallskip
	
	\noindent
	\textbf{Corollary B.} \textit{Let $R$ be a Noetherian ring containing $\mathbb{Q}$ and $A$ an $R$-subalgebra such that $A^{[n]}=A[\underline{T}] = R[\UL{X}] = R^{[n+2]}$ where $\UL{X} = (X_1, X_2, \cdots, X_{n+2})$ and $\UL{T} = (T_1, T_2, \cdots, T_n)$ are two sequence of indeterminates. Let $F \in A$ be such that $A\otimes_R \Qt{R/P} = \Qt{R/P}[F]^{[1]}$ for each minimal prime ideal $P$ of $R$. Then, the following are equivalent.}
	
	\begin{enumerate}
		\item [\rm (I)] $(F_{X_1},F_{X_2}, \cdots, F_{X_{n+2}})A[\UL{T}] = A[\UL{T}]$.
		\item [\rm (II)] \textit{$F$ is a residual coordinate of $A$.}
		\item [\rm (III)] $A = R[F]^{[1]} = R^{[2]}$.
		\item [\rm (IV)] \textit{$\JacDer{F, \UL{T}}{\UL{X}}$ is a $R$-LND with slice.}
		\item [\rm (V)] \textit{$\JacDer{F, \UL{T}}{\UL{X}}$ is a fixed point free $R$-LND.}
		\item [\rm (VI)] \textit{$I = A[\UL{T}]$ where $I$ is the ideal generated by $(n+1) \times (n+1)$ minors of $\JacMat{F, \UL{T}}{\UL{X}}$.}
	\end{enumerate}
	
	\medskip
	
	It is to be observed that, in Corollary B, the $R$-algebra $A$ is a retract of $R^{[n+2]}$, and therefore, by Corollary \ref{Cor_Neena-Nikhilesh-Sagnik_Retraction-A2Fib} $A$ is an $\A{2}$-fibration over $R$. 
	
	\medskip
	
	Kahoui-Ouali, in \cite{Kahoui_A2-fib_triviality-criterion}, established a criterion for a stably polynomial $\A{2}$-fibration to be a polynomial algebra (see Theorem \ref{Kahoui-Ouali-Thm}), which is closely related to Corollary B. In Section \ref{Sec_Kahoui-Ouali-Result} we revisit the result of Kahoui-Ouali in view of Corollary B. At this point we must mention that the result of Kahoui-Ouali, i.e., Theorem \ref{Kahoui-Ouali-Thm}, helped as one of the sources to build intuitions for subsequent development of the subject (see \cite{Babu-Das_Struct_A2-fib_FPF-LND}, \cite{Babu-Das-Lokhande_Rank-and-Rigidity}).
	
	\section{Preliminaries}
	
	In this section we recall definitions and quote some results required for the subsequent sections. 
	
	\begin{defn}
		Let $R \subseteq A$ be rings. $R$ is called a \textit{retract} of $A$ if there exists an $R$-algebra surjection $\phi: A \longrightarrow R$ such that $\phi(x) =x$ for each $x \in R$.
		%
		%
	\end{defn}
	
	\begin{defn}		
		Let $R$ be a ring, $A$ an $R$-algebra. and $D: A \longrightarrow A$ be an $R$-derivation. 
		$D$ is defined to be \textit{fixed point free} if $D(A)A = A$.
		$D$ is called a \textit{locally nilpotent $R$-derivation} ($R$-LND), if for each $x \in A$, there exists $n \in \mbbN$ such that $D^n(x) =0$. 	
	\end{defn}	
	
	We now quote the required results. The first one is by Hamann (\cite[Theorem 2.8]{Haman_Invariance}).
	
	\begin{thm} \label{Hamann}
		Let $R$ be a Noetherian ring containing $\mathbb{Q}$ and $A$ an $R$-algebra such that $A^{[m]} = R^{[m+1]}$ for some $m \in \mathbb{N}$. Then, $A = R^{[1]}$.
	\end{thm}
	
	The following L\"{u}roth-type result is by Abhyankar-Eakin-Heinzer (\cite[Theorem 4.1]{AEH_Coff}) and Russell-Sathaye (\cite[Corollary 3.4]{RS_FIND}).
	
	\begin{thm} \label{AEH_Luroth-type}
		If $R \subseteq A \subseteq R^{[n]}$ are UFDs such that $\trdeg{R}{A} =1$, then $A=R^{[1]}$.
	\end{thm}
	
	A criterion by Russell-Sathaye (\cite[Theorem 2.3.1]{RS_FIND}) for an algebra to be a polynomial algebra states as
	
	\begin{thm} \label{Russell-Sathaye_Criterion}
		Let $R \subseteq A$ be integral domains such that $A$ is finitely generated over $R$. Suppose, there exists $\pi \in R$ which is prime in $A$ such that $\pi A \cap R= \pi R$, $A[1/\pi] =R[1/\pi]^{[1]}$ and the image of $R/\pi R$ is algebraically closed in $A/\pi A$. Then, $A = R^{[1]}$.
	\end{thm}
	
	The next result is by Sathaye (\cite[Theorem 1]{Sat_Pol-two-var-DVR}).
	
	\begin{thm} \label{Sathaye_DVR}
		Let $R$ be a DVR containing $\mbbQ$ and $A$ an $\A{2}$-fibration over $R$. Then, $A = R^{[2]}$.
	\end{thm}

	Asanuma established the following structure theorem (\cite[Theorem 3.4]{Asanuma_fibre_ring}) of affine fibrations over Noetherian rings.
	
	\begin{thm} \label{Asanuma_struct-fib-th}
		Let $R$ be a Noetherian ring and $B$ an $\A{r}$-fibration over $R$. Then, $\Omega_R(B)$ is a projective $B$-module of rank $r$ and $B$ is an $R$-subalgebra (up to an isomorphism) of a polynomial ring $R^{[m]}$ for some $m \in \mathbb{N}$ such that $B^{[m]}=\Sym{R^{[m]}} {\Omega_R(B) \otimes_B R^{[m]}}$. Therefore, $B$ is a retract of $R^{[n]}$ for some $n$.
	\end{thm}
	
	The below result of Asanuma-Bhatwadekar describes the structure of an $\A{2}$-fibration over a one dimensional Noetherian ring containing $\mbbQ$ (\cite[Theorem 3.8 \& Remark 3.13]{Asan-Bhatw_Struct-A2-fib}).
	
	\begin{thm} \label{Asanuma-Bhatwadelar_A2-fib-Structure}
		Let $R$ be a one-dimensional Noetherian ring containing $\mbbQ$ and $A$ an $\A{2}$-fibration over $R$. Then there exists $H \in A$ such that $A$ is an $\A{1}$-fibration over $R[H]$.
	\end{thm}
	
	The following is a result by Essen (\cite{Essen_Around-Cancellation}; also see\cite{Essen-Maubach-Berson_Der-div-zero}, \cite{BD_LND}, \cite{Daigle-Freudenburg_UFD-LND-Rank-2}, \cite{Rentchler_Operations-du-Groupe})
	
	\begin{thm} \label{Essen-fpf-LND-A2}
		Let $R$ be a ring containing $\mathbb{Q}$, $A = R[X,Y]$ and $D$ a fixed point free $R$-LND of $A$. Then, $\KerD = R^{[1]}$ and $A = \KerD^{[1]}$.
	\end{thm}
	
	%
	%
	%
	
	Below we quote a result on $\A{1}$-forms by Das (see \cite[Theorem 3]{Das_A1-form}).
	
	\begin{thm} \label{Das_A1-form}
		Let $R$ be a ring and $A$ be an $R$-algebra such that
		\begin{enumerate}
			\item [\rm (i)]$A$ is a UFD.
			\item [\rm (ii)]There is a retraction $\Phi : A \longrightarrow R$.
			\item [\rm (iii)]There exists a faithfully flat ring homomorphism $\eta: R \longrightarrow R'$ such that $A \otimes_R R' = R'^{[1]}$.
		\end{enumerate}
		Then, $A  = R^{[1]}$.
	\end{thm}

	\begin{rem}
		Though, in \cite{Das_A1-form}, Theorem \ref{Das_A1-form} has been proven assuming $\eta: R \longrightarrow R'$ to be a finite and faithfully flat ring homomorphism, one can easily see that the finiteness condition is redundant.
	\end{rem}
	
	The following result on residual coordinate is by Das-Dutta (\cite[Corollary 3.6, Lemma 3.12, Theorem 3.16 \& Corollary 3.19]{DD_residual}).
	
	\begin{thm}\label{DD_fib}
		Let $R$ be a Noetherian ring and $A$ an $\A{n}$-fibration over $R$. Suppose, $\UL{W} \in A$ is an $m$-tuple residual coordinate of $A$. Then, $A$ is an $\A{n-m}$-fibration over $R[\UL{W}]$ and $\Omega_R(A) = \Omega_{R[\UL{W}]}(A) \oplus A^m$. Further, if $A$ is stably polynomial over $R \hookleftarrow \mbbQ$ and $n-m=1$, then $A = R[\UL{W}]^{[1]} = R^{[n]}$.
	\end{thm}
	
	It is to be noted that though Das-Dutta, in \cite{DD_residual}, proved Theorem \ref{DD_fib} (see \cite[Corollary 3.19]{DD_residual}) with the hypothesis that the base ring is a Noetherian domain containing $\mbbQ$, from their proof it follows that Theorem \ref{DD_fib} holds over Noetherian rings (not necessarily domains) containing $\mbbQ$.
	
	The next result is by Chakraborty-Dasgupta-Dutta-Gupta (\cite[Theorem 5.9]{Neena-Sagnik-Nikhilesh_Retract-of-polynomial}).
	
	\begin{thm} \label{Neena-Nikhilesh-Sagnik_Retraction-A2Fib}
		Let $R$ be a Noetherian domain containing $\mbbQ$ and $A$ an integral domain containing $R$ such that $\trdeg{R}{A} =2$. Then, $A$ is an $\A{2}$-fibration over $R$ if and only if  $A$ is an $R$-algebra retract of $R^{[m]}$ for some integer $m$.
	\end{thm}
	
	We end this section by observing that Theorem \ref{Neena-Nikhilesh-Sagnik_Retraction-A2Fib} extends over Noetherian rings too.
	
	\begin{cor}
		\label{Cor_Neena-Nikhilesh-Sagnik_Retraction-A2Fib}
		Let $R$ be a Noetherian ring containing $\mbbQ$ and $A$ an over ring of $R$ such that $\trdeg{R/P}{A \otimes_R R/P} =2$ for each minimal prime ideal $P$ of $R$. Then, the following statements are equivalent.
		\begin{enumerate}
			\item [\rm (I)] $A$ is an $R$-algebra retract of $R^{[m]}$ for some integer $m$.
			\item [\rm (II)] $A$ is an $\A{2}$-fibration over $R$.
		\end{enumerate}
	\end{cor}
	
	\begin{proof}
		\uline{(I) $\implies$ (II):} We assume that (I) holds. Fix a minimal prime ideal $P_0$ of $R$. Clearly, $A \otimes_R R/P_0$ is a retract of $(R/P_0)^{[m]}$ and hence $A \otimes_R R/P_0$ is finitely generated and faithfully flat over $R/P_0$. Now, by Theorem \ref{Neena-Nikhilesh-Sagnik_Retraction-A2Fib}, we see that $A \otimes_R R/P_0$ is an $\A{2}$-fibration over $R/P_0$, and therefore, $A \otimes_{R} k(\bar{Q}) = (A \otimes_{R} R/P_0) \otimes_{R/P_0} k(\bar{Q}) = k(\bar{Q})^{[2]}$ for all $\bar{Q} \in \Spec{R/P_0}$. Since any $\bar{Q} \in \Spec{R/P_0}$ is of the form $Q/P_0$ for some $Q \in \Spec{R}$, we see that $k(\bar{Q}) = \Qt{(R/P_0)/\bar{Q}} = \Qt{R/Q} = k(Q)$, and therefore, we have $A \otimes_{R} k(Q) = (R \otimes_R k(Q)^{[2]}$ for any $Q \in \Spec{R}$ such that $P_0 \subseteq Q$. Now, since $P_0$ is arbitrary, we have $A \otimes_{R} k(P) = k(P)^{[2]}$ for each $P \in \Spec{R}$. By the hypothesis $A$ is finitely generated and faithfully flat over $R$, and hence it  follows that $A$ is an $\A{2}$-fibration over $R$.
		
		\medskip
		
		\uline{(II) $\implies$ (I):} Follows directly from Theorem \ref{Asanuma_struct-fib-th}. 
	\end{proof}

	\section{Recognizing residual coordinates of $\A{2}$-fibrations}
	\label{Sec_Recognizing-residual-variables}
	In this section we prove Theorem A and Corollary B. First, we prove a special case of Theorem A.
	
	\begin{prop} \label{Prop_Recognize-ResCord-of-A2Fib}
		Let $R$ be a Noetherian domain containing $\mathbb{Q}$ with quotient field $K$ and $A$ an $R$-algebra such that $A$ is a retraction of $B = R[X_1, X_2, \cdots, X_{n}] = R^{[n]}$ and $\trdeg{R}{A} =2$. Let $F \in A$ be such that $A\otimes_R K = K[F]^{[1]}$. If $(F_{X_1},F_{X_2}, \cdots, F_{X_{n}})B = B$, then $F$ is a residual coordinate of $A$.
	\end{prop}
	
	\begin{proof} Note that, by Theorem \ref{Neena-Nikhilesh-Sagnik_Retraction-A2Fib}, $A$ is an $\A{2}$-fibration over $R$. Since an element $h \in A$ is a residual coordinate of $A$ over $R$ if and only if $h$ is a residual coordinate of $A_P$ over $R_P$ for all $P \in \Spec{R}$, without loss of generality we assume $R$ to be local with maximal ideal $\m$ and residue field $k = R/\m$. Since $A \otimes_R K = K[F]^{[1]}$, we only need to show that $A \otimes_R k(P) = k(P)[\ol{F}]^{[1]}$ for each $P \in \Spec{R} \backslash \{0\}$ where $\ol{F}$ denotes the image of $F$ in $A \otimes_R k(P)$. We prove this using induction on $d:= \dim(R)$.
		
		\smallskip
		
		If $d = 0$, there is nothing to prove as $F$ is already a coordinate of $A \otimes_R K = K^{[2]}$. Assume that $d = 1$. Since $R$ is a one-dimensional Noetherian local domain, by the Krull-Akizuki theorem there exists a discrete valuation ring (DVR) $(C,\pi)$ such that $R \subseteq C \subseteq K$ and the residue field $L = C/({\pi})$ is finite over $k$. By Theorem \ref{Sathaye_DVR} we get $A \otimes_R C = C^{[2]}$. We shall show that $A \otimes_R C = C[F]^{[1]}$. Since $F$ is a generic coordinate in $A \otimes_R  C$, by Theorem \ref{Russell-Sathaye_Criterion} it is enough to show that $L[\ol{F}]$ is algebraically closed in $A \otimes_R L = (A \otimes_R C) \otimes_C L = L^{[2]}$. By Theorem \ref{AEH_Luroth-type} the algebraic closure of $L[\ol{F}]$ in $A \otimes_R L = L^{[2]}$ is of the form $L[U]$. So, we have $L[\ol{F}] \subseteq L[U] \subseteq A \otimes_R L \subseteq B \otimes_R L = L[X_1, X_2, \cdots, X_{n}]$, and therefore, the partial derivatives of $U$ with respect to $X_1, X_2, \cdots, X_{n}$ are well defined. Let us write $\ol{F}$ = $a_0 + a_1U + a_2U^2 +  \cdots  + a_mU^m$ where $a_i \in L$ and in that case we get the following equations
		$$
		\frac{\partial F}{\partial X_i} = a_1\frac{\partial U}{\partial X_i} + 2Ua_2\frac{\partial U}{\partial X_i} +  \cdots  + mU^{m-1}a_m\frac{\partial U}{\partial X_i}		
		$$
		where $i  = 1,2, \cdots , n+2$. Since $(F_{X_1},F_{X_2}, \cdots, F_{X_{n}})B = B$, there exists $b_1,b_2, \cdots ,b_{n+2}$ such that $b_1\frac{\partial F}{\partial X_1} + b_2\frac{\partial F}{\partial X_2} + \cdots + b_{n}\frac{\partial F}{\partial X_{n}} = 1$ which along with the above equation gives us the following
		$$(b_1\frac{\partial U}{\partial X_1} + b_2\frac{\partial U}{\partial X_2} +  \cdots  + b_{n}\frac{\partial U}{\partial X_{n}})(a_1 + 2a_2U +  \cdots  + ma_mU^{m-1}) = 1$$
		
		On comparing the degrees of $U$ in the above equation, we get $a_i = 0$  for all $i = 2,3, \cdots ,m$. This shows that $\ol{F}$ is linear in $U$, and therefore, $L[\ol{F}] = L[U]$. So, we have $A \otimes_R C = C[F]^{[1]}$, and hence $A \otimes_R L = L[\ol{F}]^{[1]}$. Since $L$ is a finite separable extension of $k$, we get $A \otimes_R k = k[\ol{F}]^{[1]}$, i.e., $F$ is a residual coordinate of $A$.
		
		\smallskip
		
		Now, assume that the result holds for all domains of dimension $d \le m-1$. Let $\dim(R) =m$. Consider $R_1:= R/P$ where $P$ is an height-one prime ideal of $R$ and set $A_1 := A \otimes_R R_1$. Clearly, $A_1$ is an $\mathbb{A}^2$-fibration over $R_1$. Since $\dim(R_1) \le m - 1$, by induction hypothesis $\ol{F}$ is a residual coordinate of $A_1$, and therefore, $A_1 \otimes_{R_1} k(\ol{Q}) = k(\ol{Q})[F]^{[1]}$ for all $\ol{Q} \in \Spec{R_1}$ where $\ol{Q}$ denote image of $Q \in \Spec{R}$ in $R_1$ such that $P \subseteq Q$. Since $A_1 = A \otimes_R R_1$ and $k(\ol{Q}) = k(Q)$ we have $A \otimes_R k(Q) = k(Q)[\ol{F}]^{[1]}$ for all $Q \in \Spec{R}$ such that $P \subseteq Q$. This shows we have $A \otimes_R k(P) = k(P)[\ol{F}]^{[1]}$ for all $P \in \Spec{R}$ with $\Ht{P} \ge 1$. Since $F$ is already a generic coordinate of $R[X,Y]$, i.e., $R[X,Y] \otimes_{R} k(0) = k(0)[F]^{[1]}$, it follows that $F$ is a residual coordinate of $A$.
	\end{proof}
	
	We now prove Theorem A.
	
	\begin{thm} \label{Thm_Recognize-ResCord-of-A2Fib-RingCase}
		Let $R$ be a Noetherian ring containing $\mathbb{Q}$ and $A$ an $R$-algebra such that $A$ is a retraction of $B = R[X_1, X_2, \cdots, X_{n}] = R^{[n]}$ and $F$ an element of $A$. Suppose, $\trdeg{R/P}{A \otimes_{R} R/P} =2$ and $A\otimes_R \Qt{R/P} = \Qt{R/P}[F]^{[1]}$ for each for each minimal prime ideal $P$ of $R$. If $(F_{X_1},F_{X_2}, \cdots, F_{X_{n}})B = B$, then $F$ is a residual coordinate of $A$.
	\end{thm}
	
	\begin{proof}
		Fix a minimal prime ideal $P_0$ of $R$. Clearly, $R/P_0$ is Noetherian and $A \otimes_R R/P_0$ is a retract of $B \otimes_R R/P_0$, and therefore, by Proposition \ref{Prop_Recognize-ResCord-of-A2Fib} we see that $F$ is a residual coordinate of $A\otimes_R R/P_0$ over $R/P_0$, i.e., $(A \otimes_{R} R/P_0) \otimes_{R/P_0} k(\bar{Q}) = ((R/P_0)[F] \otimes_{R/P_0} k(\bar{Q}))^{[1]}$ for all $\bar{Q} \in \Spec{R/P_0}$. Since any $\bar{Q} \in \Spec{R/P_0}$ is of the form $Q/P_0$ for some $Q \in \Spec{R}$ such that $P_0 \subseteq Q$, we see that $k(\bar{Q}) = \Qt{(R/P_0)/\bar{Q}} = \Qt{R/Q} = k(Q)$, and therefore, we have $A \otimes_{R} k(Q) = (R[F] \otimes_R k(Q))^{[1]}$ for any $Q \in \Spec{R}$ such that $P_0 \subseteq Q$. Now, since $P_0$ is arbitrary, we have $A \otimes_{R} k(P) = (R[F] \otimes_{R} k(P))^{[1]}$ for each $P \in \Spec{R}$, i.e., $F$ is a residual coordinate of $A$.
	\end{proof}
	
	Next, we prove Corollary B.
	\begin{cor} \label{Cor_Recognize-ResCord-of-StablyTrivA2Fib}
Let $R$ be a Noetherian ring containing $\mathbb{Q}$ and $A$ an $R$-subalgebra such that $A^{[n]}=A[\underline{T}] = R[\UL{X}] = R^{[n+2]}$ where $\UL{X} = (X_1, X_2, \cdots, X_{n+2})$ and $\UL{T} = (T_1, T_2, \cdots, T_n)$ are two sequence of indeterminates. Let $F \in A$ be such that $A\otimes_R \Qt{R/P} = \Qt{R/P}[F]^{[1]}$ for each minimal prime ideal $P$ of $R$. Then, the following are equivalent.

\begin{enumerate}
	\item [\rm (I)] $(F_{X_1},F_{X_2}, \cdots, F_{X_{n+2}})A[\UL{T}] = A[\UL{T}]$.
	\item [\rm (II)] $F$ is a residual coordinate of $A$.
	\item [\rm (III)] $A = R[F]^{[1]} = R^{[2]}$.
	\item [\rm (IV)] $\JacDer{F, \UL{T}}{\UL{X}}$ is a $R$-LND with slice.
	\item [\rm (V)] $\JacDer{F, \UL{T}}{\UL{X}}$ is a fixed point free $R$-LND.
	\item [\rm (VI)] $I = A[\UL{T}]$ where $I$ is the ideal of $A[\UL{T}]$ generated by the $(n+1) \times (n+1)$ minors of $\JacMat{F, \UL{T}}{\UL{X}}$.
\end{enumerate}
\end{cor}

\begin{proof} \UL{(I) $\implies$ (II):} Assume that (I) holds. Clearly, $A$ is a retract of $R[\UL{X}]$ and $\trdeg{R/P}{A \otimes_{R} R/P} =2$ for each minimal prime ideal $P$ of $R$. Now, the implication follows from Theorem \ref{Thm_Recognize-ResCord-of-A2Fib-RingCase}.\\
\UL{(II) $\implies$ (III):} Directly follows from Theorem \ref{DD_fib}.\\
\UL{(III) $\implies$ (IV):} Assume $A = R[F, G] = R^{[2]}$ for some $G \in A$. Then, $A [\UL{T}] = R[F, G , \UL{T}] = R[\UL{X}] = R^{[n+2]}$ and hence $\det(\JacMat{F, G, \UL{T}}{\UL{X}}) \in R^*$. Now one can see that  $\Delta: = \JacDer{F, \UL{T}}{\UL{X}} = \det(\JacMat{F, G, \UL{T}}{\UL{X}}) \JacDer{F, \UL{T}}{F, G, \UL{T}}$ is an $R$-LND and $\Delta(G) \in R^*$, i.e., $\Delta$ has a slice.\\
\UL{(IV) $\implies$ (V):} Trivial.\\
\UL{(V) $\implies$ (VI):}  Assume that $\Delta: = \JacDer{F, \UL{T}}{\UL{X}}$ is fixed point free, i.e., $(\Delta(X_1), \Delta(X_2), \cdots, \Delta(X_{n+2}))A[\UL{T}] = A[\UL{T}]$. Now, the implication follows from the equality $I = (\Delta(X_1), \Delta(X_2), \cdots, \Delta(X_{n+2}))A[\UL{T}]$ which one can check very easily.\\
\UL{(VI) $\implies$ (I):} Assume that (VI) holds. Clearly, $I \subseteq (F_{X_1},F_{X_2}, \cdots, F_{X_{n+2}})A[\UL{T}]$ and hence the implication follows.
\end{proof}

The following two observations are consequences of Theorem \ref{Thm_Recognize-ResCord-of-A2Fib-RingCase}. These results may help determine residual coordinates of affine $n$-fibrations.
	
	\begin{cor} \label{Cor_Triviality-of-An-fib_Partial-der_general}
	Let $R$ be a Noetherian ring containing $\mathbb{Q}$ and $A$ an $\mathbb{A}^n$-fibration over $R$. Suppose that $\UL{W} = (W_1, W_2, \cdots, W_{n-2})$ is a $(n-2)$-tuple residual coordinate of $A$ and $F \in B$ be such that $A \otimes_R \Qt{R/P} = (R[W_1, W_2, \cdots, W_{n-2}, F] \otimes_R \Qt{R/P})^{[1]}$ for each minimal prime ideal $P$ of $R$. Then the following hold.

\begin{enumerate}
	\item [\rm (I)] $A$ is an $\A{2}$-fibration over $R[\UL{W}]$ and there exist $\UL{X} = (X_1, X_2, \cdots, X_{m})$, $m \ge 2$ such that $A$ is a retract of $R^{[m]}$.
	
	\item [\rm (II)] If $(F_{X_1},F_{X_2}, \cdots, F_{X_{m}})A[\UL{T}] = A[\UL{T}]$, then $(W_1, W_2, \cdots, W_{n-2}, F)$ is a residual coordinate tuple of $A$.
\end{enumerate}	
	\end{cor}
	
\begin{proof}
	\UL{(I):} Since $\UL{W}$ is a residual coordinate tuple of $A$, by Theorem \ref{DD_fib} it directly follows that $A$ is an $\A{2}$-fibration over $R[\UL{W}]$, and hence by Corollary \ref{Cor_Neena-Nikhilesh-Sagnik_Retraction-A2Fib} there exist coordinate tuples $\UL{X} = (X_1, X_2, \cdots, X_{m})$, $m \ge 2$ such that $A$ is a retract of $R[\UL{X}]$.
	
	\medskip
	
	\UL{(II):} We assume that $(F_{X_1},F_{X_2}, \cdots, F_{X_{m}})A[\UL{T}] = A[\UL{T}]$. Fix a minimal prime ideal $P$ of $R$. Since by hypothesis $A \otimes_R \Qt{R/P} = (R[W_1, W_2, \cdots, W_{n-2}, F] \otimes_R \Qt{R/P})^{[1]}$ we have $A \otimes_{R[\UL{W}]} \Qt{\frac{R}{P}[\UL{W}]} = \Qt{\frac{R}{P}[\UL{W}]}[F]^{[1]}$. It is easy to see that the  minimal prime ideal of $R[X]$ are exactly of the form $P[X]$ where $P$ is a minimal prime of $R$, and therefore, by Theorem \ref{Thm_Recognize-ResCord-of-A2Fib-RingCase} it follows that $F$ is a residual coordinate of $A$ over $R[\UL{W}]$. Let $P \in \Spec{R}$. Since $R[\UL{W}]$ is faithfully flat over $R$, there exists $Q \in \Spec{R[\UL{W}]}$ such that $P = Q \cap R$. Clearly, $k(Q)$ is a field extension of $k(P)$. Note $A \otimes_{R} k(P) = k(P)^{[n]}$ is a UFD. Now, since $A \otimes_{R[\UL{W}]} k(Q) = k(Q)[F]^{[1]}$, we see that $(A \otimes_R k(P)) \otimes_{k(P)} k(Q) = ((R[F, \UL{W}] \otimes_R k(P))\otimes_{k(P)} k(Q))^{[1]}$. Therefore, by Theorem \ref{Das_A1-form} we have $A \otimes_R k(P) = k(P) [F, \UL{W}]^{[1]}$, and hence it follows that $(W_1, W_2, \cdots, W_{n-2}, F)$ is a residual coordinate tuple of $A$.
\end{proof}

	\begin{cor} \label{Cor_Triviality-of-An-fib_Partial-der_special}
	Let $R$ be a Noetherian domain containing $\mathbb{Q}$ and $A$ an $R$-algebra such that $A^{[m]} = A[\UL{T}] = R[\UL{X}] = R^{[m+n]}$ where $\UL{X} = (X_1, X_2, \cdots, X_{m+n})$ and $\UL{T} = (T_1, T_2, \cdots, T_{m})$. Suppose that $\UL{W} = (W_1, W_2, \cdots, W_{n-2})$ is a $(n-2)$-tuple residual coordinate of $A$ and $F \in B$ be such that $A \otimes_R \Qt{R/P} = (R[W_1, W_2, \cdots, W_{n-2}, F] \otimes_R \Qt{R/P})^{[1]}$ for each minimal prime ideal $P$ of $R$. Then the following are equivalent.

\begin{enumerate}
	\item [\rm (I)] $(F_{X_1},F_{X_2}, \cdots, F_{X_{n+m}})A[\UL{T}] = A[\UL{T}]$.
	\item [\rm (II)] $(W_1, W_2, \cdots, W_{n-2}, F)$ is a residual coordinate tuple of $A$.
	\item [\rm (III)] $A = R[W_1, W_2, \cdots, W_{n-2}, F]^{[1]} = R^{[n]}$.
\end{enumerate}
	\end{cor}

\begin{proof}
	Follows from Corollary \ref{Cor_Triviality-of-An-fib_Partial-der_general} and Corollary \ref{Cor_Recognize-ResCord-of-StablyTrivA2Fib}.
\end{proof}
	
	\section{An alternative proof to Kahoui-Ouali's theorem: Stably polynomial $\A{2}$-fibrations having fixed point free LNDs are polynomial algebras} \label{Sec_Kahoui-Ouali-Result}
	In this section, in view of Corollary \ref{Cor_Recognize-ResCord-of-StablyTrivA2Fib}, we revisit the main result of Kahoui-Ouali in \cite{Kahoui_A2-fib_triviality-criterion}, which states as follows (\cite[Theorem 3.1]{Kahoui_A2-fib_triviality-criterion}).
	
	\begin{thm} \label{Kahoui-Ouali-Thm}
		Let $R$ be a ring containing $\mathbb{Q}$, $A$ a stably polynomial $\A{2}$-fibration over $R$ and $D$ an $R$-LND of $A$. Then, $D$ is fixed point free if and only if $A = R^{[2]}$.
	\end{thm}
	
	Before going to the proof of Theorem \ref{Kahoui-Ouali-Thm} we first observe some results.
	\begin{lem} \label{Lem_ideal-equality}
		Let $R \subset A$ be rings, $I$ an ideal of $R$ and $J_1 \subseteq J_2$ ideals of $A$ such that $I \subset J_i \bigcap R$ for $i = 1,2$. Suppose that images of $J_1$ and $J_2$ are the same in $A/IA$, then $J_1 = J_2$.
	\end{lem}
	\begin{proof} For any element $x \in A$, let the notation $\bar{x}$ denotes the image of $x$ in $A \otimes_R R/I$. Let $x \in J_2$. Then $\bar{x} \in J_2/ IA = J_1/IA$. This shows that there exists $y \in J_1$ such that $\bar{x} - \bar{y} = \bar{0}$ in $A/IA$, i.e., $x-y \in IA$, i.e., $x \in J_1 + IA$. Since $IA \subseteq J_1$, it is clear that $x \in J_1$, and therefore we have  $J_1 = J_2$.
	\end{proof}
	
	%
	\begin{lem} \label{Lem_stably-pol-is-pol}
		Let $R \subseteq A$ be rings and $\eta$ be the nilradical of $R$. Suppose that $A^{[m]} = R^{[m+n]}$ and $A/\eta A= (R/\eta)^{[n]}$. Then, $A = R^{[n]}$.
	\end{lem}
	
	\begin{proof} For the proof, one may see \cite[Lemma 4.2]{Kahoui_A2-fib_triviality-criterion}. However, for the convenience of the reader here we sketch it. Let $X_1, X_2, \cdots, X_{m+n}$ and $T_1, T_2, \cdots, T_{m}$ be such that $A^{[m]} = A[T_1, T_2, \cdots, T_{m}] = R[X_1, X_2, \cdots, X_{m+n}] = R^{[m+n]}$. Since $A/\eta A = (R/\eta)^{[n]}$, there exist $x_1, x_2, \cdots, x_n \in A$ such that $A = R[x_1, x_2, \cdots, x_n] + \eta A$, and therefore,  $A = R[x_1, x_2, \cdots, x_n] + \eta^{\ell} A$ for all $\ell \in \mathbb{N}$. Since $\eta$ is nilpotent, one can see that $A = R[x_1, x_2, \cdots, x_n]$. Now, if $x_1, x_2, \cdots, x_n \in A$ are algebraically dependent over $R$, they remain the same in $A[T_1, T_2, \cdots, T_{m}] = A^{[m]}$, which is a contradiction to the fact that \\$A[T_1, T_2, \cdots, T_{m}]=  R[x_1, x_2, \cdots, x_n, T_1, T_2, \cdots, T_{m}] = R[X_1, X_2, \cdots, X_{m+n}] = R^{[m+n]}$.
	\end{proof}
	
	\begin{lem} \label{Lem_A2-fib-over-0dim-is-trivial}
		Let $R$ be a reduced Artinian ring and $A$ an $\A{n}$-fibration over $R$. Then, $A = R^{[n]}$.
	\end{lem}
	\begin{proof}
		Since $R$ is a zero dimensional reduced Noetherian ring, it follows that $R = \bigoplus_{i =1}^{r} k_i$ where $k_i = R/\mathfrak{m}_i$ and $\{ \mathfrak{m}_i :  i = 1,2, \cdots, r\}$ is the set of maximal ideals of $R$. Since $A$ is an $\A{n}$-fibration over $R$, for each $i = 1,2 \cdots, r$, we have $A \otimes_{R} k(\mathfrak{m}_i) = A \otimes_{R} k_i = R[X_{1i}, X_{2i}, \cdots, X_{ni}] \otimes_{R} k_i = k_i^{[n]}$ for some $X_{1i}, X_{2i}, \cdots, X_{ni} \in A \backslash R$. By prime avoidance lemma, we can choose $p_i \in \mathfrak{m}_i$ such that $p_i \notin (\bigcup_{j =1}^r \mathfrak{m}_j) \backslash \mathfrak{m}_i$. Let $q_i = (\prod_{j=1}^{r} p_j) / p_i$. For each $j = 1,2, \cdots, n$, set $X_j: = \sum_{i=1}^{r}q_i X_{ji}$. Clearly, $R[X_1, X_2, \cdots, X_n] \otimes_{R} k_i = R[X_{1i}, X_{2i}, \cdots, X_{ni}] \otimes_{R} k_i = k_i^{[n]}$ for all $i =1,2,\cdots, r$. observe that $A = A \otimes_{R} \bigoplus_{i =1}^{r} k_i = \bigoplus_{i =1}^{r} (R[X_1, X_2, \cdots, X_n] \otimes_{R} k_i) = \bigoplus_{i =1}^{r} k_i [X_1, X_2, \cdots, X_n] = (\bigoplus_{i =1}^{n} k_i )[X_1, X_2, \cdots, X_n] = R[X_1, X_2, \cdots, X_n]$. We shall show that $R[X_1, X_2, \cdots, X_n] = R^{[n]}$. Now, since $R$ is reduced and $R[X_1, X_2, \cdots, X_n] \otimes_{R} k_i = k_i^{[n]}$ for each $i = 1,2 \cdots, r$, one can see that $X_1, X_2, \cdots, X_n$ can not be algebraically dependent over $R$, and hence $A = R^{[n]}$.
	\end{proof}
	
	\begin{prop} \label{Prop_Strut-LND_FPF-or-grade2}
		Let $R$ be a Noetherian ring containing $\mathbb{Q}$ with total quotient ring $K$, $A$ an $R$-algebra such that $A^{[n]} = A[\UL{T}] = R[\UL{X}] = R^{[n+2]}$ where $\UL{X} = (X_1, X_2, \cdots, X_{n+2})$ and $\UL{T} = (T_1, T_2, \cdots, T_n)$ and $D: A \longrightarrow A$ a fixed point free $R$-LND. Then, there exists an irreducible $F \in A$ such that $A \otimes_{R} K = K[F]^{[1]} = K^{[2]}$, $\Ker{D} \otimes_R K = K[F]$ and $(F_{X_1},\cdots, F_{X_{n+2}})A[\UL{T}] = A [\UL{T}]$.
	\end{prop}
	
	\begin{proof}
		Since $A^{[n]} = A[\UL{T}] = R[\UL{X}] = R^{[n+2]}$, by Corollary \ref{Cor_Neena-Nikhilesh-Sagnik_Retraction-A2Fib} we see that $A$ is an $\A{2}$-fibration over $R$. We shall use induction on $\dim(R)$ to prove our claim. Suppose that $\dim(R) = 0$. We shall show that $A = R^{[2]}$. In view of Lemma \ref{Lem_stably-pol-is-pol} we can assume that $R$ is reduced, and hence by Lemma \ref{Lem_A2-fib-over-0dim-is-trivial} we have $A = R^{[2]}$. Since $D$ is fixed point free, by Theorem \ref{Essen-fpf-LND-A2} we get $A = \Ker{D}[s] = \Ker{D}^{[1]}$. Now, by Theorem \ref{Hamann} we have $\KerD = R[F]$ for some $F \in A$. So, we have $R^{[n+2]} = A[\UL{T}] = R[F, \UL{T}][s]$ and therefore, it follows that $(F_{X_1},\cdots, F_{X_{n+2}})A[\UL{T}] = A [\UL{T}]$. Now, assume that $\dim(R) =1$. Then, by Theorem \ref{Asanuma-Bhatwadelar_A2-fib-Structure} there exists $H$ such that $A$ is an $\A{1}$-fibration over $R[H]$. This shows that $H$ is a residual coordinate of $A$ and therefore, by Theorem \ref{DD_fib} we see that $A = R[H]^{[1]} = R^{[2]}$. Now, repeating the previous arguments we have an irreducible $F \in A$ such that $(F_{X_1},\cdots, F_{X_{n+2}})A[\UL{T}] = A [\UL{T}]$.
		
		\medskip
		
		Next, we assume that $\dim(R) = \ell$ and the result holds for any ring of dimension at most $\ell -1$. Note that there exits $F, G \in A$ such that $A[\UL{T}] \otimes_R K = K[F,G][\UL{T}]= K [\UL{X}]$ and $\Ker{D} \otimes_R K = K[F]$, and therefore, $(F_{X_1},\cdots, F_{X_{n+2}})A[\UL{T}] \otimes_R K = A[\UL{T}] \otimes_R K$. Since $R$ is Noetherian, we can choose both $F$ and $G$ to be irreducible. Now, since the ideal $(F_{X_1},\cdots, F_{X_{n+2}})A[\UL{T}]$ is finitely generated, there exists a non-zero divisor $t \in R$ such that $(F_{X_1},\cdots, F_{X_{n+2}})A[\UL{T}] [1/t] = A[\UL{T}] \otimes_R R[1/t]$, and therefore, there exits a smallest non-negative integer $m \in \mathbb{N} \cup \{ 0 \}$ such that $a_1F_{X_1} + a_2 F_{X_2} +  \cdots + a_{n+2} F_{X_{n+2}} = t^m$ where $t \in R$ and $a_i \in A[\UL{T}]$ for all $i = 1, 2, \cdots, n+2$. Set $\displaystyle I : = R \cap  (F_{X_1}, F_{X_2}, \cdots,F_{X_{n+2}})A[\UL{T}]$. Clearly, $t^m \in I$. Note that $I = R$ if and only if $(F_{X_1},\cdots, F_{X_{n+2}})A[\UL{T}] = A [\UL{T}]$ if and only if $m =0$. We claim that $(F_{X_1},\cdots, F_{X_{n+2}})A[\UL{T}] = A [\UL{T}]$.
		
		\medskip
		
		If possible let $(F_{X_1},\cdots, F_{X_{n+2}})A[\UL{T}] \ne A [\UL{T}]$ for any such choice of $F$ made above. Then, $m \ne 0$, and hence $t^m \in R$ is a non-unit and non-zero divisor. Note that $\dim(R/t^mR) = \ell -1$, and hence by the induction hypothesis $(F_{X_1}, \cdots, F_{X_{n+2}} )A[\UL{T}] \otimes_R R/t^mR = A[\UL{T}] \otimes_{R}R/t^mR$. Now, since $t^m \in (F_{X_1}, \cdots, F_{X_{n+2}} )A[\UL{T}] \subseteq A[\UL{T}]$, by Lemma \ref{Lem_ideal-equality} we have $(F_{X_1}, \cdots, F_{X_{n+2}} )A[\UL{T}] = A[\UL{T}]$ which is a contradiction.
	\end{proof}
	
	We now give an proof to Theorem \ref{Kahoui-Ouali-Thm}. 
	
	\medskip
	
	\textbf{Proof to Theorem \ref{Kahoui-Ouali-Thm}:} Since a stably polynomial algebra $A$ over a ring $R$ is a polynomial algebra over $R$ if and only if $A \otimes_R R/\eta$ is a polynomial algebra over $R/\eta$ where $\eta$ is the nilradical of $R$ (see Lemma \ref{Lem_stably-pol-is-pol}), without loss of generality we may assume that $R$ is reduced. By another reduction technique (see \cite[Lemma 4.3 ]{Kahoui_A2-fib_triviality-criterion}) we can further assume that $R$ is Noetherian. Now, the consecutive applications of Proposition \ref{Prop_Strut-LND_FPF-or-grade2}, Theorem \ref{Cor_Recognize-ResCord-of-StablyTrivA2Fib} and Theorem \ref{DD_fib} establish the result.
	
	\section{Example} \label{Sec_Example}
	
	In this section we discuss two well known examples to demonstrate the applications of Corollary \ref{Cor_Recognize-ResCord-of-StablyTrivA2Fib} and Proposition \ref{Prop_Recognize-ResCord-of-A2Fib} (or Theorem \ref{Thm_Recognize-ResCord-of-A2Fib-RingCase}). The first one is by Hochster (\cite{Hochster_nonunique-coeff}).
	
	\begin{ex} \label{Ex_Hochster}
		Let $R = \mathbb{R}[X,Y,Z]/(X^2 + Y^2 + Z^2 -1)$ and  $K$ denote the quotient field of $R$. Let $x, y \text{ and } z$ denote the image of the coordinates $X,Y \text{ and }Z$ respectively in $R$. Let $s:= xU + yV + zW \in  R[U,V,W] = R^{[3]}$. It is easy to see that $R[U,V,W] = A[s] = A^{[1]}$ where $A = R[U - xs, \ V - ys, \ W - zs]$. Define a map $\phi: R[U,V,W] \longrightarrow A$ by $(U, \ V, \ W) \mapsto (U - xs, \ V - ys, \ W - zs)$. One can check that $\phi$ is an retraction onto $A$ and $\Ker{\phi} = sR[U,V,W]$. By Theorem \ref{Neena-Nikhilesh-Sagnik_Retraction-A2Fib} we see that $A$ is an $\A{2}$-fibration over $R$. Set $g := zV - yW = z(V - ys) - y(W- zs)\in A$ and $h := (y U - x V) = y(U - xs) - x(V - ys) \in A$. One can check that $A \otimes_{R} K = K[g,h]$, i.e., $g$ and $h$ form a generic coordinate pair of $A$. However, we see that $(\frac{\partial}{\partial U}(g), \frac{\partial}{\partial V}(g), \frac{\partial}{\partial W}(g)) \ne R[U,V,W]$ and $(\frac{\partial}{\partial U}(h), \frac{\partial}{\partial V}(h), \frac{\partial}{\partial W}(h)) \ne R[U,V,W]$, and therefore, by Corollary \ref{Cor_Recognize-ResCord-of-StablyTrivA2Fib} neither $g$ nor $h$ are residual coordinate of $A$. 
	\end{ex}
	
	Next, we consider an example of an $\A{2}$-fibration by Asanuma-Bhatwadekar (\cite[Example 3.12]{Asan-Bhatw_Struct-A2-fib}).
	
	\begin{ex} \label{Ex_Asanuma-Bhatwadekar}
		Set $S := \mathbb{C}[t,X]( = \mathbb{C}^{[2]})$, $R_1 := \mathbb{C}[t^2, t^3]$ and $R := R_1[X]( = R_1^{[1]})$. Let $A = R[Y + tX^2Y^2]+
		(t^2, t^3)R[Y] \subseteq S[Y]$. By a result of Greither (\cite[Theorem 3.2]{Greither_Seminormal-projective-invertible}) it follows that $A$ is a retract of $R[Y_1, Y_2, \cdots, Y_n] = R^{[n]}$ for some $n \in \mathbb{N}$. Now, since $R = R_1[X]$, it follows that $A$ is a retract of $R_1[X, Y_1, Y_2, \cdots, Y_n] = R^{[n+1]}$ and $\trdeg{R_1}{A} =2$, and therefore, by Theorem \ref{Neena-Nikhilesh-Sagnik_Retraction-A2Fib} we see that $A$ is an $\A{2}$-fibration over $R_1$. Clearly, $X \in A$ and $(\frac{\partial}{\partial X}(X), (\frac{\partial}{\partial Y_1}(X), \cdots, (\frac{\partial}{\partial Y_n}(X)) = R[X, Y_1, Y_2, \cdots, Y_n]$, and hence by Proposition \ref{Prop_Recognize-ResCord-of-A2Fib} it follows that $X$ is a residual coordinate of $A$.
		
	\end{ex}
	\section*{Acknowledgment} 
	The second author acknowledges SERB, Govt. of India for their MATRICS grant under the file number MTR/2022/000247.

\bibliographystyle{alpha}
	\normalem
	\bibliography{reference}

\begin{thebibliography}{CDDG21}

\bibitem[AB97]{Asan-Bhatw_Struct-A2-fib}
Teruo Asanuma and S.~M. Bhatwadekar.
\newblock Structure of {${\bf A}\sp 2$}-fibrations over one-dimensional
  {N}oetherian domains.
\newblock {\em J. Pure Appl. Algebra}, 115(1):1--13, 1997.

\bibitem[AEH72]{AEH_Coff}
Shreeram~S. Abhyankar, Paul Eakin, and William Heinzer.
\newblock On the uniqueness of the coefficient ring in a polynomial ring.
\newblock {\em J. Algebra}, 23:310--342, 1972.

\bibitem[Asa87]{Asanuma_fibre_ring}
Teruo Asanuma.
\newblock Polynomial fibre rings of algebras over {N}oetherian rings.
\newblock {\em Invent. Math.}, 87(1):101--127, 1987.

\bibitem[BD93]{BD_RES}
S.~M. Bhatwadekar and Amartya~K. Dutta.
\newblock On residual variables and stably polynomial algebras.
\newblock {\em Comm. Algebra}, 21(2):635--645, 1993.

\bibitem[BD97]{BD_LND}
S.~M. Bhatwadekar and Amartya~K. Dutta.
\newblock Kernel of locally nilpotent {$R$}-derivations of {$R[X,Y]$}.
\newblock {\em Trans. Amer. Math. Soc.}, 349(8):3303--3319, 1997.

\bibitem[BD21]{Babu-Das_Struct_A2-fib_FPF-LND}
Janaki~Raman Babu and Prosenjit Das.
\newblock Structure of {$\Bbb A^2$}-fibrations having fixed point free locally
  nilpotent derivations.
\newblock {\em J. Pure Appl. Algebra}, 225(12):Paper No. 106763, 12, 2021.

\bibitem[BDL21]{Babu-Das-Lokhande_Rank-and-Rigidity}
Janaki~Raman Babu, Prosenjit Das, and Swapnil~A. Lokhande.
\newblock Rank and rigidity of locally nilpotent derivations of affine
  fibrations.
\newblock {\em Comm. Algebra}, 49(12):5214--5228, 2021.

\bibitem[BvM01]{Essen-Maubach-Berson_Der-div-zero}
Joost {Berson}, Arno {van den Essen}, and Stefan {Maubach}.
\newblock {Derivations having divergence zero on $R[X,Y]$.}
\newblock {\em {Isr. J. Math.}}, 124:115--124, 2001.

\bibitem[CDDG21]{Neena-Sagnik-Nikhilesh_Retract-of-polynomial}
Sagnik Chakraborty, Nikhilesh Dasgupta, Amartya~Kumar Dutta, and Neena Gupta.
\newblock Some results on retracts of polynomial rings.
\newblock {\em Journal of Algebra}, 567:243--268, 2021.

\bibitem[Das12]{Das_A1-form}
Prosenjit Das.
\newblock A note on factorial {$\Bbb A^1$}-forms with retractions.
\newblock {\em Comm. Algebra}, 40(9):3221--3223, 2012.

\bibitem[Das15]{Das_cancel}
Prosenjit Das.
\newblock On cancellation of variables of the form $bt^n-a$ over affine normal
  domains.
\newblock {\em J. Pure Appl. Algebra}, 219(12):5280--5288, 2015.

\bibitem[DD14]{DD_residual}
Prosenjit Das and Amartya~K. Dutta.
\newblock A note on residual variables of an affine fibration.
\newblock {\em J. Pure Appl. Algebra}, 218(10):1792--1799, 2014.

\bibitem[DF98]{Daigle-Freudenburg_UFD-LND-Rank-2}
Daniel Daigle and Gene Freudenburg.
\newblock Locally nilpotent derivations over a \text{UFD} and an application to
  rank two locally nilpotent derivations of $k[x_1, \cdots , x_n]$.
\newblock {\em J. Algebra}, 204(2):353--371, 1998.

\bibitem[EK13]{Kahoui_Residual}
M'hammed El~Kahoui.
\newblock On residual coordinates and stable coordinates of {$R^{[3]}$}.
\newblock {\em Arch. Math. (Basel)}, 100(1):35--41, 2013.

\bibitem[EKEO21]{Kahoui-Essamaoui-Ouali_ResCoord-one-dim}
M'hammed El~Kahoui, Najoua Essamaoui, and Mustapha Ouali.
\newblock Residual coordinates over one-dimensional rings.
\newblock {\em J. Pure Appl. Algebra}, 225(6):Paper No. 106629, 10, 2021.

\bibitem[EKO16]{Kahoui_A2-fib_triviality-criterion}
M'hammed El~Kahoui and Mustapha {Ouali}.
\newblock {A triviality criterion for $\mathbb{A}^2$-fibrations over a ring
  containing $\mathbb{Q}$.}
\newblock {\em {J. Algebra}}, 459:272--279, 2016.

\bibitem[Gre81]{Greither_Seminormal-projective-invertible}
Seminormality, projective algebras, and invertible algebras.
\newblock {\em Journal of Algebra}, 70(2):316--338, 1981.

\bibitem[Ham75]{Haman_Invariance}
Eloise Hamann.
\newblock On the {$R$}-invariance of {$R[X]$}.
\newblock {\em J. Algebra}, 35:1--16, 1975.

\bibitem[{Hoc}72]{Hochster_nonunique-coeff}
M.~{Hochster}.
\newblock {Nonuniqueness of coefficient rings in a polynomial ring.}
\newblock {\em {Proc. Am. Math. Soc.}}, 34:81--82, 1972.

\bibitem[Lah19]{Lahiri_partial-coordinate-system}
Animesh Lahiri.
\newblock A note on partial coordinate system in a polynomial ring.
\newblock {\em Comm. Algebra}, 47(3):1099--1101, 2019.

\bibitem[Ren68]{Rentchler_Operations-du-Groupe}
Rudolf Rentschler.
\newblock Op\'{e}rations du groupe additif sur le plan affine.
\newblock {\em C. R. Acad. Sci. Paris S\'{e}r. A-B}, 267:A384--A387, 1968.

\bibitem[RS79]{RS_FIND}
Peter Russell and Avinash Sathaye.
\newblock On finding and cancelling variables in {$k[X,\,Y,\,Z]$}.
\newblock {\em J. Algebra}, 57(1):151--166, 1979.

\bibitem[Sat83]{Sat_Pol-two-var-DVR}
Avinash Sathaye.
\newblock Polynomial ring in two variables over a {DVR}: a criterion.
\newblock {\em Invent. Math.}, 74(1):159--168, 1983.

\bibitem[vdE07]{Essen_Around-Cancellation}
Arno van~den Essen.
\newblock Around the cancellation problem.
\newblock In {\em Affine algebraic geometry}, pages 463--481. Osaka Univ.
  Press, Osaka, 2007.

\end{thebibliography}
\end{document}